# A fully new path to prove Riemann Hypothesis

**J. M. Zhu**

## Abstract

In the past 160 years, people used to study Zeta function only regarding it as a complex function. Generally, complex functions are much more complex than real functions, and are hard to graph. So, people cannot grasp the nature of them easily. Therefore, it may be a promising way to try to correspond Zeta function to real function so that we can return to the real domain to study RH.

In fact, under Laplace transform, $F(s)=L[f(x)]=\int_0^\infty f(x)e^{-sx}dx$, the whole picture of Zeta function is very clear and simple, and the problem can be greatly simplified. And by Laplace transform, most integral and convolution operations can be converted into algebraic operations, which greatly simplifies calculating and analysis. The first part of this paper points out ***the essence of Zeta function*** displayed under Laplace transform, and the second part obtains ***the error function*** $Je(x)-lie(x)$ ***which equivalent to the error function*** $J(x)-li(x)$ ***got by Riemann before, while this new error function is much simpler than the Riemann's.*** The third part estimates the maximum absolute value of the new error function by two ways, thus ***proving the famous equivalent proposition of Riemann Hypothesis:***

$$|li(x)-\pi(x)|=O(x^{0.5+\varepsilon}),\quad \varepsilon>0,\ \pi(x) \text{ is the prime-counting function},\ li(x)=\int_0^x \frac{dt}{\log t},\quad t\neq 1$$

The fourth part makes some other meaningful discussions and obtains other valuable results.

## 1. The essence of Zeta function

In the original definition, Zeta function is $\zeta(s)=1+\frac{1}{2^s}+\frac{1}{3^s}+\frac{1}{4^s}+\dots,\quad \mathrm{Re}(s)>1$

Riemann extended the Zeta function to the region where s≠1 using analytical extension. New Zeta function is in the form of contour integration, which appears simple but is actually more inconvenient to analyze than the original Zeta function. ***The original Zeta function already contains all the information about the distribution of prime numbers. All Zeta functions mentioned later are in original definition and Re(s)>1.*** If we use **Laplace inverse transform** to view Zeta function, a clear picture is shown:

$$\zeta(s)=1+\frac{1}{2^s}+\frac{1}{3^s}+\frac{1}{4^s}+\dots=1+e^{-s\log 2}+e^{-s\log 3}+e^{-s\log 4}+\dots$$

***This represents a series of Dirac delta functions at the points of x=0, log2, log3, log4, ..., which can be recorded as：*** $zeta(x)=L^{-1}[\zeta(s)]=\delta(x)+\delta(x-\log 2)+\delta(x-\log 3)+\delta(x-\log 4)+\dots$

It may be still difficult to understand what this means, but once it is integrated, the truth is clear.:

$$zeta_1(x)=\int_0^x zeta(t)dt=L^{-1}[\frac{1}{s}\zeta(s)]=u(x)+u(x-\log 2)+u(x-\log 3)+u(x-\log 4)...$$

***This represents a combination of unit step functions start at the points x=0, log2, log3, log4, ..., as shown in Figure 1, obviously has an exact upper bound function*** $f_z^+(x)=e^x$ ***and a lower bound function*** $f_z^-(x)=e^x-1$. It can be regarded as a combination of step functions with constant ordinate and logarithmic contraction from abscissa. This is why many of the series decomposition expressions of $\zeta(s)$ contain $\frac{1}{s-1}$, because the main item of real function corresponding to $\frac{\zeta(s)}{s}$ is $e^x$, $L(e^x)=\frac{1}{s-1}$.

Then we let: $r(x)=e^x-zeta_1(x)$, obviously at $x=\log(N+c)$, $0\le c<1$, $r(x)=c$.

We define $R(s)=L[r(x)]=\frac{1}{s-1}-\frac{1}{s}\zeta(s)$, **...... (1.1)**

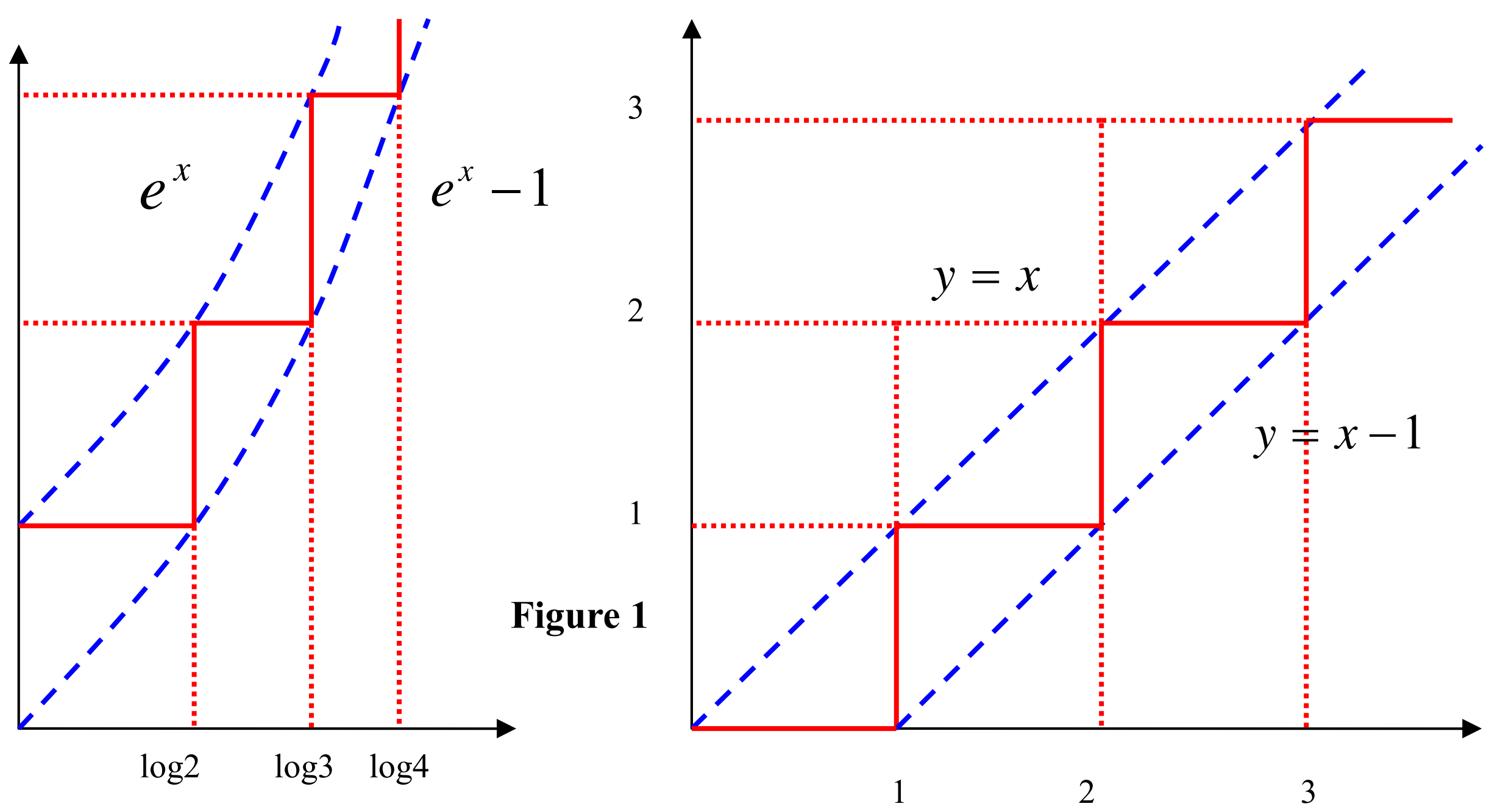

**Figure 1**

And it is easy to find: $L(eta(x))=\eta(s)=(1-\frac{2}{2^s})\zeta(s)=1-\frac{1}{2^s}+\frac{1}{3^s}-\frac{1}{4^s}+\ldots$, $\int_0^x eta(t)dt$ as shown in Figure 2, is the contraction of the horizontal coordinates of the unit square wave function:

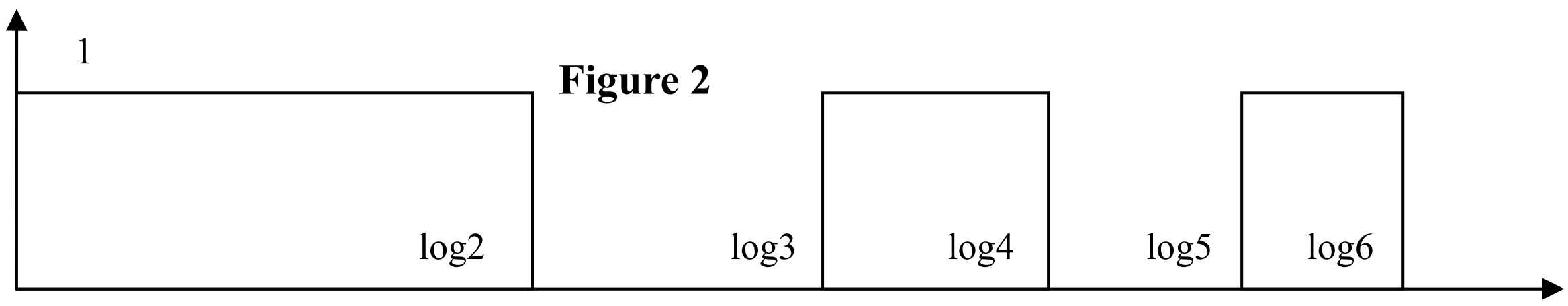

**Figure 2**

## 2. Get the error function of $Je(x) - lie(x)$

In his famous paper, Riemann used $J(x)$ function to calculate $\pi(x)$. The definition of $J(x)$ is:

$$J(x) = \sum_p \frac{1}{k} u(x - p^k) = u(x-2) + u(x-3) + \frac{1}{2}u(x-2^2) + u(x-5) + u(x-7) + \frac{1}{3}u(x-2^3)$$

$+\frac{1}{2}u(x-3^2) + u(x-11) + u(x-13) + \frac{1}{4}u(x-2^4) + \ldots$, where $p = prime\ number$.

So, by Mobius inverse transform, there is:

$$\pi(x) = J(x) - \frac{1}{2}J(x^{\frac{1}{2}}) - \frac{1}{3}J(x^{\frac{1}{3}}) - \frac{1}{5}J(x^{\frac{1}{5}}) + \frac{1}{6}J(x^{\frac{1}{6}}) - \frac{1}{7}J(x^{\frac{1}{7}}) + \ldots \qquad \textbf{...... (2.1)}$$

Obviously, $J(x) - \pi(x) \approx \frac{1}{2}J(x^{\frac{1}{2}})$ is the main error of $J(x) - \pi(x)$. In other words, if we can prove $|J(x) - li(x)| = O(x^{0.5+\varepsilon})$, it is equivalent to prove $|li(x) - \pi(x)| = O(x^{0.5+\varepsilon}), \quad \varepsilon > 0$.

We define ***new functions:*** $Je(x) = J(e^x)$, $Je(\log x) = J(x)$, and $JE(s) = L[Je(x)]$.

$Je(x)$ can be derived from $\zeta(s)$. By Euler product formula, $\zeta(s) = \Pi(1 - p^{-s})^{-1}$, the two sides are logarithmic and multiplied by $\frac{1}{s}$, $\frac{1}{s}Je(s) = \frac{1}{s}\log\zeta(s) = -\frac{1}{s}\log\Pi(1-p^{-s}) = \frac{1}{s}(\frac{1}{2^s} + \frac{1}{2\cdot 2^s}$

$+\frac{1}{3\cdot 2^{3s}} + \ldots + \frac{1}{3^s} + \frac{1}{2\cdot 3^{2s}} + \frac{1}{3\cdot 3^{3s}} + \ldots + \frac{1}{5^s} + \frac{1}{2\cdot 5^{2s}} + \frac{1}{3\cdot 5^{3s}} + \ldots), \quad Je(x) = \sum_p \frac{1}{k}u(x - \log p^k)$.

It is the same as $J(x)$'s constant ordinate and logarithmic contraction from abscissa.

For comparing with $li(x)$, we must make the same transform on $li(x)$. $li(x)$ has the series decomposition expression as following, where $\gamma = 0.5772\ldots$ is Euler constant:

$li(x) = \gamma + \log|\log x| + \log x + \frac{\log^2 x}{2\cdot 2!} + \frac{\log^3 x}{3\cdot 3!} + \frac{\log^4 x}{4\cdot 4!} + \ldots$ The logarithmic contraction of the abscissa of $li(x)$ is to change the $\log x$ of the upper form to $x$, or $lie(x) = li(e^x)$:

$lie(x) = \gamma + \log x + x + \frac{x^2}{2\cdot 2!} + \frac{x^3}{3\cdot 3!} + \frac{x^4}{4\cdot 4!} + \ldots$ Make Laplace transform to $lie(x)$:

$$LE(s) = L(lie(x)) = \frac{\gamma}{s} - \frac{1}{s}(\gamma + \log s) + \frac{1}{s^2} + \frac{1}{2s^3} + \frac{1}{3s^4} + \frac{1}{4s^5} + \ldots = -\frac{1}{s}\log s - \frac{1}{s}\log(1 - \frac{1}{s})$$

$= -\frac{1}{s}\log(s-1)$. And by: $\int_1^x \frac{e^t}{t}dt = \int_1^x (1+t+\frac{t^2}{2!}+\frac{t^3}{3!}+...)\frac{dt}{t} = \log x + x + \frac{x^2}{2\cdot 2!}+\frac{x^3}{3\cdot 3!}+...$

$-(1+\frac{1}{2\cdot 2!}+\frac{1}{3\cdot 3!}+...) \approx lie(x)$, $lie(x)$ is the integral of $\frac{e^x}{x}$.

If $|J(x)-li(x)| = O(x^{0.5+\varepsilon_1}), \varepsilon_1 > 0$, that requires $|Je(x)-lie(x)| = O(e^{0.5x+\varepsilon_2}), \varepsilon_2 > 0$

Error function: $L[Je(x)-lie(x)] = \frac{1}{s}\log\zeta(s) + \frac{1}{s}\log(s-1) = \frac{1}{s}\log((s-1)\zeta(s))$

$= \frac{1}{s}\log\left((s-1)\frac{\zeta(s)}{s}\right) + \frac{1}{s}\log s$. For $L^{-1}\left(\frac{1}{s}\log s\right) = -\gamma - \log x$, the order is very low, if the main item of error function is $O\left(e^{0.5x+\varepsilon_2}\right), \varepsilon_2 > 0$, $-\log x - \gamma$ can be ignored. So we pay attention at:

$$ER(s) = \frac{1}{s}\log\left((s-1)\frac{\zeta(s)}{s}\right) = \frac{1}{s}\log\left((s-1)(\frac{1}{s-1} - R(s))\right) = \frac{1}{s}\log(1-(s-1)R(s)) = -\frac{1}{s}(s-1)R(s)$$

$$-\frac{1}{2s}(s-1)^2R^2(s) - \frac{1}{3s}(s-1)^3R^3(s) - \frac{1}{4s}(s-1)^4R^4(s) - ... - \frac{1}{ks}(s-1)^kR^k(s) - ... \quad \textbf{...... (2.2)}$$

## 3. Proof of Riemann Hypothesis

Although we have got the **(2.2)** formula, it is still very difficult to inverse transform it directly, no mention to other research. Because the function represented by this exact expression contains many vibration and irregular components, so it is too hard to simplify. ***Fortunately, Riemann Hypothesis does not require the exact expression of the error function, only the good upper and lower bounds of*** $\pi(x)$ ***are needed, which are both simple functions and contain no any vibration and irregular components. It provides the possibility for some equivalent methods. Therefore, the standard upper and lower bound functions of*** $Je(x) - lie(x)$ ***also do not contain vibration and irregular terms, so they can be approximated by some simple functions.***

Next I will provide ***two ways to estimate the maximum absolute value of the error function.*** But at first, I will explain more on the reasons why I do so.

i) We have 2 methods to calculate $(A+B)^2$: (1) $(A+B)^2 = A\cdot A + A\cdot B + B\cdot A + B\cdot B$, calculate each item and then sum them up; (2) calculate $A+B=C$, and then $(A+B)^2 = C\cdot C$. In this case, if we do as method (1), that is very complex, even impossible. So method (2) is the only choice.

ii) What we need to estimate is the maximum absolute value of each $L^{-1}[\frac{1}{ks}(s-1)^kR^k(s)]$.

If a function $G(s)$ can satisfies $\left|L^{-1}\left(\frac{1}{ks}(s-1)^k R^k(s)\right)\right|_{\max} \le L^{-1}\left(\frac{1}{ks}G^k(s)\right)$, then we call $G(s)$ as a ***convolution equivalent function***. Although $G(s) \ne (s-1)R(s)$, for our target, it is useful and enough.

**Proof 1**: $L^{-1}[\frac{1}{ks}(s-1)^k R^k(s)] = L^{-1}[\frac{s^k}{ks}\cdot\frac{(s-1)^k R^k(s)}{s^k}] < L^{-1}[\frac{s^k}{ks}\cdot\frac{G(s)}{s^k}]$, we try to get $\frac{G(s)}{s}$ from $\frac{1}{s}(s-1)R(s) = (1-\frac{1}{s})R(s) = -\left(\frac{R(s)}{s} - R(s)\right)$. $\frac{R(s)}{s}$ is equivalent to the integral of $r(x)$, as in Figure 3. Since $x = \log N$, $y = N$, so the area $x \cdot y = N \log N$. Then we look at the rectangles: the bottom D2 length is log2, height is 1; then the upper one D3 length is log3, height is 1; ..., and the top Dn length is logN, height is 1. So the " ‖ " shadow area at $x = \log N$ is (the error $O(1/N^3)$ is ignored):

$$S = x \cdot y - (D2 + D3 + D4 + .. + Dn) = N \log N - \sum_{m=1}^{N} \log m,$$

Using the formula: $\sum_{m=1}^{N} \log m = N \log N - N + \frac{\log N}{2} + \frac{\log 2\pi}{2} + \frac{1}{12N} + O\left(\frac{1}{N^3}\right)$.

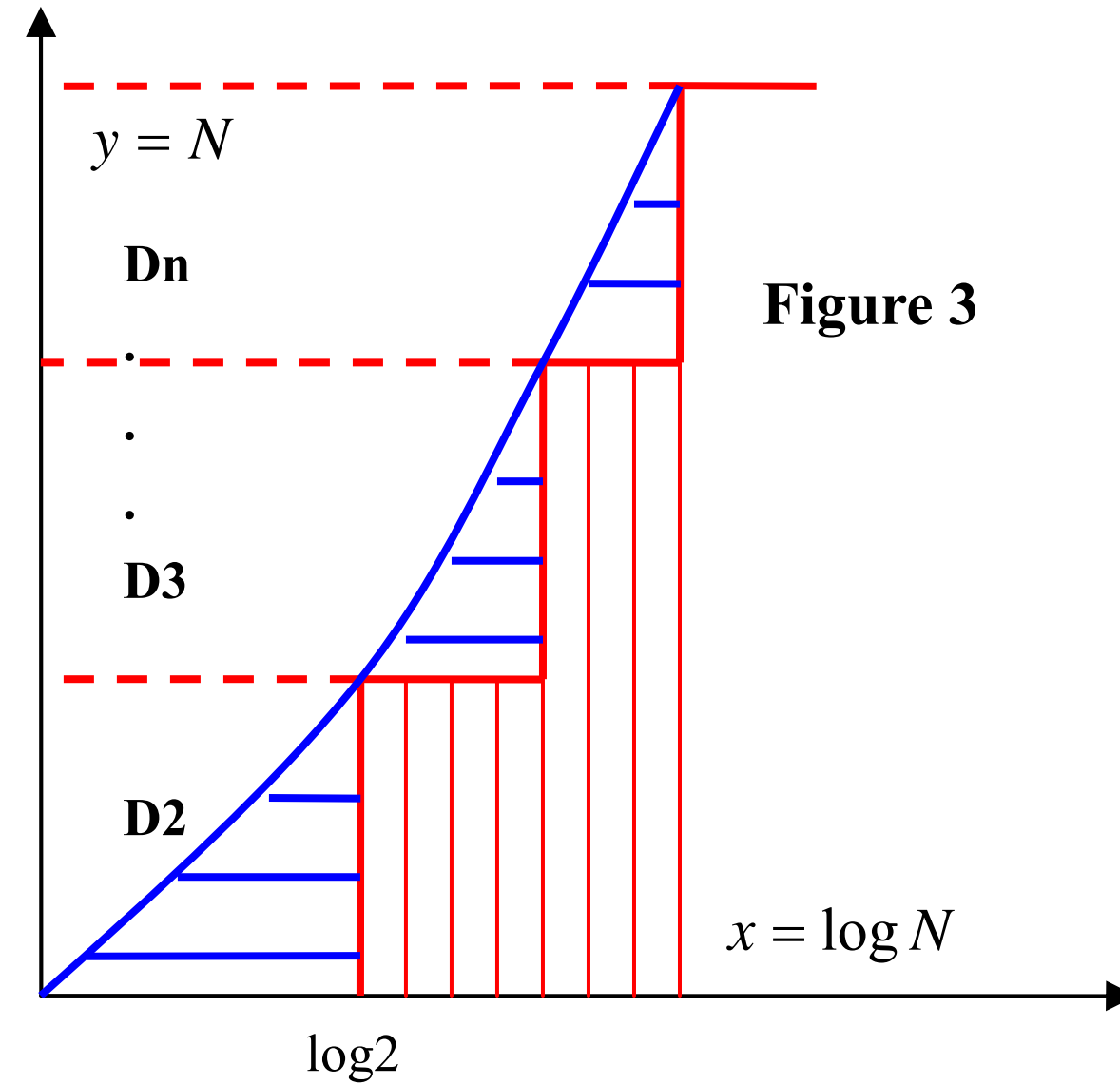

**Figure 3**

$$S = N \log N - \sum_{m=1}^{N} \log m = N \log N - (N \log N - N + \frac{\log N}{2} + \frac{\log 2\pi}{2} + \frac{1}{12N})$$

$$= N - \frac{\log N}{2} - \frac{\log 2\pi}{2} - \frac{1}{12N}$$

The integral of $r(x)$ is the "=" shadow area:

$$\int_0^{\log N} r(x)dx = \int_0^{\log N} e^x dx - S = N - 1 - N + \frac{\log N}{2} + \frac{\log 2\pi}{2} + \frac{1}{12N} = \frac{\log N}{2} + \frac{\log 2\pi}{2} - 1 + \frac{1}{12N}$$

So at $x = \log N$, $\int_0^{\log N} r(x)dx - r(\log N) = \frac{\log N}{2} + \frac{\log 2\pi}{2} - 1 + \frac{1}{12N} - 0$, that means:

$$\int_0^x r(t)dt - r(x) = \frac{x}{2} + \frac{\log 2\pi}{2} - 1 + \frac{1}{12e^x}, \quad \left(\frac{R(s)}{s} - R(s)\right) = \frac{1}{2s^2} + \frac{\log 2\pi}{2s} - \frac{1}{s} + \frac{1}{12(s+1)}$$

$$=\frac{1}{s}\left(\frac{1}{2s}+\frac{\log 2\pi}{2}-1+\frac{s}{12(s+1)}\right)=\frac{1}{s}\left(\frac{1}{2s}+\frac{\log 2\pi}{2}-1+\frac{1}{12}-\frac{1}{12(s+1)}\right),\quad \frac{\log 2\pi}{2}-1+\frac{1}{12}=0.00227...$$

The small constant 0.00227, its effect is far less than the item $-\frac{1}{12(s+1)}$ in the case of multiple convolution, so the maximum absolute value of the item $L^{-1}[\frac{1}{ks}(s-1)^k R^k(s)]$ is:

$$L^{-1}\left|\frac{1}{ks}(s-1)^k R^k(s)\right|_{\max} \approx L^{-1}\left(\frac{1}{ks}(\frac{1}{2s}-\frac{1}{12(s+1)})^k\right) < L^{-1}\left(\frac{1}{ks}(\frac{1}{2s})^k\right), \textbf{...... (3.1)}$$

At $x=\log(N+c)$, $r(x)$ adds $c$, the integral of $r(x)$ adds very small value, $\int_0^x r(t)dt-r(x)$, the negative term is larger. The (3.1) formula is more established. We can calculate as below:

$$\int_0^{\log(N+c)} r(x)dx - r(\log(N+c)) = N+c-1-\left(S+N(\log(N+c)-\log N)\right)-c$$

$$=N-1-(N-\frac{\log N}{2}-\frac{\log 2\pi}{2}-\frac{1}{12N})-N\log(1+c/N)=\frac{\log N}{2}+\frac{\log 2\pi}{2}+\frac{1}{12N}-1-N(\frac{c}{N}-\frac{c^2}{2N^2})$$

$$=\frac{\log(N+c)}{2}-\frac{\log(1+c/N)}{2}+\frac{\log 2\pi}{2}-1-c+\frac{1}{12N}+\frac{c^2}{2N}=\frac{\log(N+c)}{2}+\frac{\log 2\pi}{2}-1-c+\frac{1-6c+6c^2}{12(N+c)}$$

(the error $O(1/N^2)$ is ignored). That means: $\int_0^x r(t)dt-r(x)=\frac{x}{2}+\frac{\log 2\pi}{2}-1-c+\frac{1-6c+6c^2}{12e^x}$

$$\frac{R(s)}{s}-R(s)=\frac{1}{2s^2}+\frac{\log 2\pi}{2s}-\frac{1+c}{s}+\frac{1-6c+6c^2}{12(s+1)}=\frac{1}{s}\left(\frac{1}{2s}+\frac{\log 2\pi}{2}-1-c+\frac{1-6c+6c^2}{12}-\frac{1-6c+6c^2}{12(s+1)}\right)$$

For $1-6c+6c^2=6(c-0.5)^2-0.5$, when $c=0.5$, $\frac{R(s)}{s}-R(s)=\frac{1}{s}\left(\frac{1}{2s}+\frac{\log 2\pi}{2}-\frac{37}{24}+\frac{1}{24(s+1)}\right)$

$=\frac{1}{s}\left(\frac{1}{2s}-1.14258+\frac{1}{24(s+1)}\right)$; when $c\to 1, c<1$, $\frac{R(s)}{s}-R(s)=\frac{1}{s}\left(\frac{1}{2s}-1.51758-\frac{1}{12(s+1)}\right)$.

We can see ***however the $c$ is changed from 0 to 1, the (3.1) formula is always established. So the maximum absolute value of the error function is:***

$$\left|L^{-1}(ER(s))\right|_{\max} < L^{-1}[\frac{1}{s}(\frac{1}{2s})+\frac{1}{2s}(\frac{1}{2s})^2+\frac{1}{3s}(\frac{1}{2s})^3+...+\frac{1}{ks}(\frac{1}{2s})^k+...]=\frac{x}{2}+\frac{(x/2)^2}{2\cdot 2!}+$$

$$\frac{(x/2)^3}{3\cdot 3!}+...+\frac{(x/2)^k}{k\cdot k!}+...<\frac{ae^{0.5x}}{x}=a\left(\frac{1}{x}+\frac{1}{2}+\frac{x}{2^2\cdot 2!}+\frac{x^2}{2^3\cdot 3!}+...+\frac{x^k}{2^{k+1}\cdot(k+1)!}+...\right)$$

Then we have $\frac{(x/2)^k}{k\cdot k!}<\frac{ax^k}{2^{k+1}\cdot(k+1)!}$, $a>\frac{2(k+1)}{k}=2+\frac{2}{k}$, let $a=3$, only the item $\frac{x}{2}>\frac{3x}{2^2 2!}$.

But we can let $\frac{x}{2}+\frac{(x/2)^2}{2\cdot 2!}+\frac{(x/2)^3}{3\cdot 3!}<\frac{3x}{2^2\cdot 2!}+\frac{3x^2}{2^3\cdot 3!}+\frac{3x^3}{2^4\cdot 4!}$, $\frac{1}{2}-\frac{3}{8}<\frac{x^2}{2^7}-\frac{x^2}{16\cdot 9}$, $x>12$. That means:

$$\left|Je(x)-lie(x)\right|_{\max}=\left|L^{-1}\left(ER(s)\right)\right|_{\max}-\log x-\gamma<\frac{3e^{0.5x}}{x},\quad \textit{for all}\ \ x>12\text{ , ...... } \mathbf{(3.2)}$$

**Proof 2:** $(s-1)R(s)=(s-1)\left(\frac{1}{s-1}-\frac{\zeta(s)}{s}\right)=1-(1-\frac{1}{s})\zeta(s)=-\left(\zeta(s)-1\right)+\frac{\zeta(s)}{s}$

The Dirac delta function can be regards as distribution density function which centralize its density at one point with its integral=1. $\zeta(s)-1=\frac{1}{2^s}+\frac{1}{3^s}+\frac{1}{4^s}+...$ is a combination of Dirac delta functions. And it is easy to find:

$$L^{-1}\left[\frac{1}{s}\left(\frac{1}{2^s}+\frac{1}{2.5^s}+\frac{1}{3^s}+\frac{1}{3.5^s}+...\right)^k\cdot\left(\frac{1}{2}\right)^k\right]<L^{-1}\left[\frac{1}{s}\left(\frac{1}{2^s}+\frac{1}{3^s}+\frac{1}{4^s}+...\right)^k\right]$$ because some density of $\zeta(s)-1=\frac{1}{2^s}+\frac{1}{3^s}+\frac{1}{4^s}+...$ are delayed. On the other hand, there should be:

$$L^{-1}\left[\frac{1}{s}\left(\frac{1}{2^s}+\frac{1}{3^s}+\frac{1}{4^s}+...\right)^k\right]<L^{-1}\left[\frac{1}{s}\left(\frac{1}{1.5^s}+\frac{1}{2^s}+\frac{1}{2.5^s}+\frac{1}{3^s}+...\right)^k\cdot\left(\frac{1}{2}\right)^k\right]$$, because some of the density of $\zeta(s)-1=\frac{1}{2^s}+\frac{1}{3^s}+\frac{1}{4^s}+...$ are advanced.

According to the same logic, we can conclude that: If the density 1 of $\frac{1}{2^s}$ is scattered to the interval $(0,\log 2]$, the density 1 of the $\frac{1}{3^s}$ is scattered to the interval $(\log 2,\log 3]$, ... , such as:

$$\frac{1}{m}\left(\left(1\frac{1}{m}\right)^{-s}+\left(1\frac{2}{m}\right)^{-s}+...\left(1\frac{m}{m}\right)^{-s}\right)+\frac{1}{m}\left(\left(2\frac{1}{m}\right)^{-s}+\left(2\frac{2}{m}\right)^{-s}+...\left(2\frac{m}{m}\right)^{-s}\right)+...$$ When $m\to\infty$, the maximum value of distribution density function equal to $e^x$ and $L\left[e^x\right]=\frac{1}{s-1}$.

On the other hand, if the density 1 of the $\frac{1}{2^s}$ is scattered to the interval $[\log 2,\log 3)$, and the density 1 of $\frac{1}{3^s}$ is scattered to the interval $[\log 3,\log 4)$, ... , such as:

$$\frac{1}{m}\left(2^{-s}+\left(2\frac{1}{m}\right)^{-s}+...\left(2\frac{m-1}{m}\right)^{-s}\right)+\frac{1}{m}\left(3^{-s}+\left(3\frac{1}{m}\right)^{-s}+...\left(3\frac{m-1}{m}\right)^{-s}\right)+...$$ When $m\to\infty$, the function takes the minimum value $e^x\cdot u(x-\log 2)$, and $L\left[e^x\cdot u(x-\log 2)\right]=\frac{2}{2^s}\frac{1}{s-1}$.

In convolution sense, the limit value has only two cases: $-\frac{1}{2^s}-\frac{1}{3^s}-\frac{1}{4^s}-...+\frac{1}{s}\zeta(s)=-\frac{1}{s-1}+\frac{1}{s-1}-R(s)=-R(s)$; or $-\frac{1}{2^s}-\frac{1}{3^s}-\frac{1}{4^s}-...+\frac{1}{s}\zeta(s)=-\frac{2}{2^s}\frac{1}{s-1}+\frac{1}{s-1}-R(s)$.

Since $\left(1-\frac{2}{2^s}\right)\frac{1}{s-1}$ represents only a small segment of function $e^x$ in the interval $[0,\log 2]$, it is obvious contrary to the symbol represented by $-R(s)$, so the absolute upper bound can only appear in the first case, that is, $\left|L^{-1}\left(\frac{1}{ks}(s-1)^k R^k(s)\right)\right|_{\max} \le L^{-1}\left(\frac{1}{ks}R^k(s)\right)$.

It is easy to prove $L^{-1}[R^k(s)] < L^{-1}\left[\frac{1}{(2s)^k}\right]$: $\frac{1}{s-1}\cdot\frac{\zeta(s)}{s} = \left(\frac{1}{s-1}-\frac{1}{s}\right)\zeta(s) = \frac{1}{s-1}\cdot\left(\frac{1}{s-1}-R(s)\right)$

$L^{-1}\left[\left(\frac{1}{s-1}-\frac{1}{s}\right)\zeta(s)\right]_{x=\log N} = N(1+\frac{1}{2}+\frac{1}{3}+...\frac{1}{N})-N = N(\log N+\gamma+\frac{1}{2N}-\frac{1}{12N^2})-N = N\log N$ $-(1-\gamma)N+\frac{1}{2}-\frac{1}{12N}$, that means: $\left(\frac{1}{s-1}-\frac{1}{s}\right)\zeta(s) = \frac{1}{(s-1)^2}-\frac{1-\gamma}{s-1}+\frac{1}{2s}-\frac{1}{12(s+1)} = \frac{1}{(s-1)^2}-\frac{R(s)}{s-1}$,

$\frac{R(s)}{s-1} = \frac{1-\gamma}{s-1}-\frac{1}{2s}+\frac{1}{12(s+1)}$, $R(s) = \frac{1}{2s}-\frac{1}{6(s+1)}+\frac{7}{12}-\gamma = \frac{1}{2s}-\frac{1}{6(s+1)}+0.0061...$ Similarly, the small constant 0.0061..., its effect under multiple convolutions is much less than $-\frac{1}{6(s+1)}$. What we need to estimate is the maximum absolute value of $L^{-1}\left(\frac{1}{s}R^k(s)\right)$. Let $F(s)=\frac{1}{2s}-\frac{1}{6(s+1)}$, $a=0.0061...$, obviously there is: $\frac{1}{2}-\frac{1}{6}=\frac{1}{3}\le L^{-1}(F(s)) = \frac{1}{2}-\frac{1}{6e^x}<\frac{1}{2}$, then we continues:

$$\frac{1}{s}R^k(s) = \frac{1}{s}(F(s)+a)^k = \frac{1}{s}\left(F^k(s)+C_k^1 aF^{k-1}(s)+C_k^2 a^2 F^{k-2}(s)+...+C_k^{k-1}a^{k-1}F(s)+a^k\right)$$

We can find its highest order is $L^{-1}\left(\frac{1}{s}F^k(s)\right) < L^{-1}\left(\frac{1}{s}\left(\frac{1}{2s}\right)^k\right)$, and 0.0061 is much less than 1, that means even $x$ is a small value, the items $\frac{1}{s}\left(C_k^1 aF^{k-1}(s)+C_k^2 a^2 F^{k-2}(s)+...+C_k^{k-1}a^{k-1}F(s)+a^k\right)$ represent the real function which much less than the highest order $L^{-1}\left(\frac{1}{s}F^k(s)\right)$.

We can conclude there is always $L^{-1}\left(\frac{1}{s}R^k(s)\right) < L\left(\frac{1}{s}\left(\frac{1}{2s}\right)^k\right)$ by other method, see the Figure 4.

The **blue line** is $r(x)$, it can be decomposed into many different sizes of rectangles. Then we get the integral of $r(x)$ is the sum of the area of all these ***small rectangles***.

$$\int_0^x r(t)dt = \frac{1}{2}(\log 2 - \log 1.5 + \log 3 - \log 2.5 + \log 4 - \log 3.5 + \ldots) + \frac{1}{4}(\log 1\frac{1}{2} - \log 1\frac{1}{4} + \log 2 - \log 1\frac{3}{4}$$

$$+ \log 2\frac{1}{2} - \log 2\frac{1}{4} + \log 3 - \log 2\frac{3}{4} + \log 3\frac{1}{2} - \log 3\frac{1}{4} + \log 4 - \log 3\frac{3}{4} + \ldots) + \ldots + \frac{1}{2^m}(\log 1\frac{2}{2^m} - \log 1\frac{1}{2^m}$$

$+ \log 1\frac{4}{2^m} - \log 1\frac{3}{2^m} + \ldots) + \ldots$ ***Note: the above and below are summation of finite items, until the item is not greater than*** $x$**.** Since there is always:

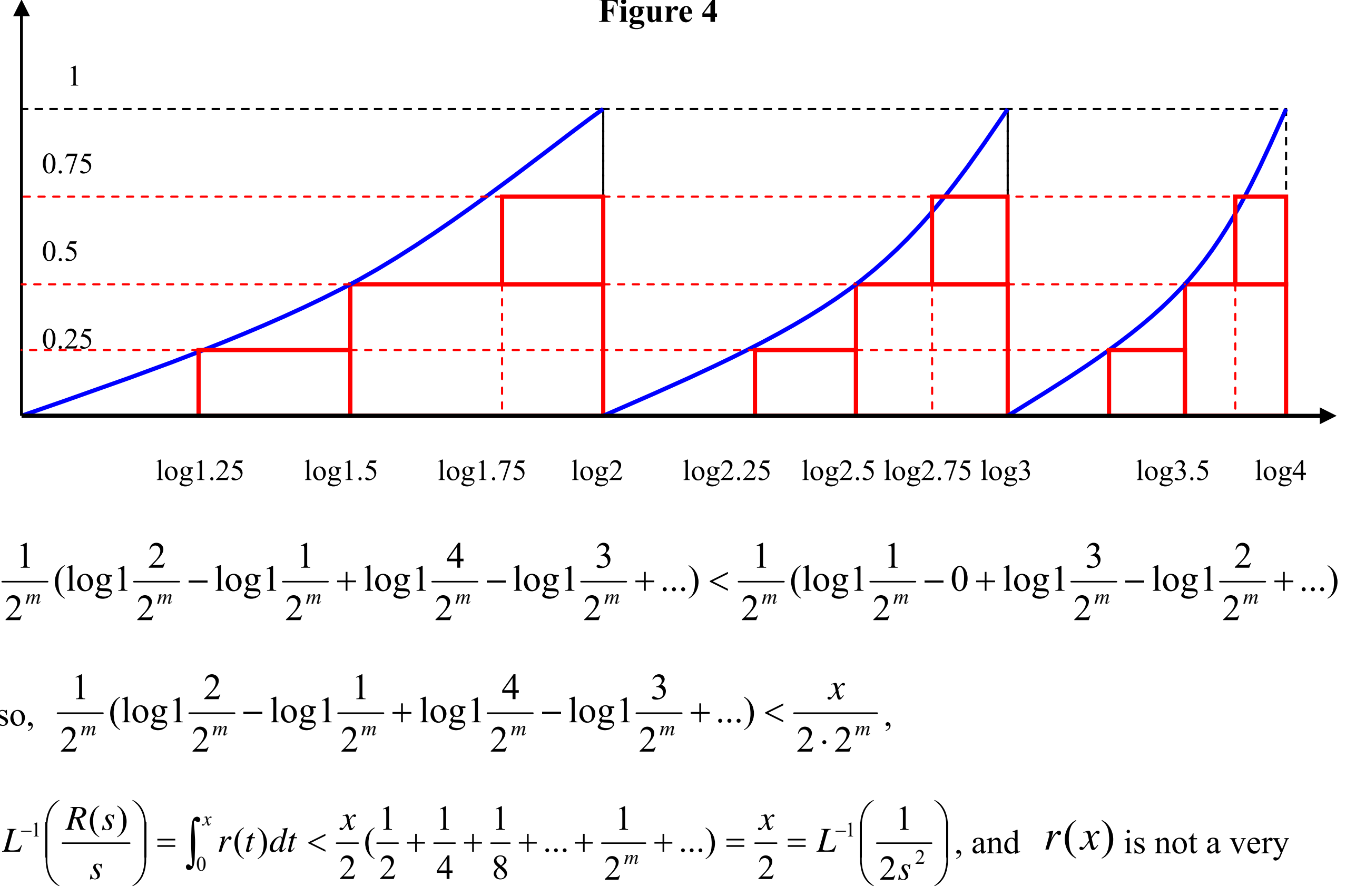

$$\frac{1}{2^m}(\log 1\frac{2}{2^m} - \log 1\frac{1}{2^m} + \log 1\frac{4}{2^m} - \log 1\frac{3}{2^m} + \ldots) < \frac{1}{2^m}(\log 1\frac{1}{2^m} - 0 + \log 1\frac{3}{2^m} - \log 1\frac{2}{2^m} + \ldots)$$

so, $\frac{1}{2^m}(\log 1\frac{2}{2^m} - \log 1\frac{1}{2^m} + \log 1\frac{4}{2^m} - \log 1\frac{3}{2^m} + \ldots) < \frac{x}{2\cdot 2^m}$,

$L^{-1}\left(\frac{R(s)}{s}\right) = \int_0^x r(t)dt < \frac{x}{2}(\frac{1}{2} + \frac{1}{4} + \frac{1}{8} + \ldots + \frac{1}{2^m} + \ldots) = \frac{x}{2} = L^{-1}\left(\frac{1}{2s^2}\right)$, and $r(x)$ is not a very vibrant function, that means the convolution equivalent real function of $r(x)$ always less than $\frac{1}{2}$, and $L^{-1}\left(\frac{1}{s}R^k(s)\right) < L^{-1}\left(\frac{1}{s}\left(\frac{1}{2s}\right)^k\right)$ always be established. The same as in **Proof 1** page 7, we can obtain:

$$\left|Je(x) - lie(x)\right|_{\max} = \left|L^{-1}\left(ER(s)\right)\right|_{\max} - \log x - \gamma < \frac{3e^{0.5x}}{x}, \quad \textit{for all} \;\; x > 12 \text{, ...... } \mathbf{(3.2)}$$

Formula **(3.2)** indicates $\left|J(x) - li(x)\right| < \frac{3\sqrt{x}}{\log x}$, from formula **(2.1)** on page 3, we know:

$J(x) - \pi(x) \approx \frac{1}{2}J(\sqrt{x}) \approx \frac{1}{2}li\left(\sqrt{x}\right)$, so we get $\left|\pi(x) + \frac{1}{2}li(\sqrt{x}) - li(x)\right| < \frac{3\sqrt{x}}{\log x}$. Compare with page7,

$$li(\sqrt{x})=lie(\log\sqrt{x})=\int_1^{\log\sqrt{x}}\frac{e^t}{t}dt=\int_1^{\log\sqrt{x}}(\frac{1}{t}+1+\frac{t}{2!}+\frac{t^2}{3!}+...)dt=\log t+t+\frac{t^2}{2\cdot 2!}+\frac{t^3}{3\cdot 3!}+...$$

$$+\frac{t^k}{k\cdot k!}+...\Big|_1^{\log\sqrt{x}},\quad \frac{b\cdot e^t}{t}=b\left(\frac{1}{t}+1+\frac{t}{2!}+\frac{t^2}{3!}+...\frac{t^k}{(k+1)!}+...\right),\quad \text{if }\ \frac{t^k}{k\cdot k!}<\frac{b\cdot t^k}{(k+1)!},\quad b>\frac{k+1}{k}=1+\frac{1}{k}.$$

So we can conclude $\dfrac{e^{\log\sqrt{x}}}{\log\sqrt{x}}=\dfrac{2\sqrt{x}}{\log x}<li(\sqrt{x})=lie(\log\sqrt{x})<\dfrac{2e^{\log\sqrt{x}}}{\log\sqrt{x}}=\dfrac{4\sqrt{x}}{\log x}$.

$$\left|\pi(x)-li(x)+\frac{1}{2}li(\sqrt{x})\right|<\frac{3\sqrt{x}}{\log x},\quad -\frac{3\sqrt{x}}{\log x}-\frac{1}{2}li(\sqrt{x})<\pi(x)-li(x)<\frac{3\sqrt{x}}{\log x}-\frac{1}{2}li(\sqrt{x})$$

$$-\frac{3\sqrt{x}}{\log x}-\frac{2\sqrt{x}}{\log x}<\pi(x)-li(x)<\frac{3\sqrt{x}}{\log x}-\frac{\sqrt{x}}{\log x},\quad -\frac{5\sqrt{x}}{\log x}<\pi(x)-li(x)<\frac{2\sqrt{x}}{\log x}.$$

**The famous equivalent proposition $|\pi(x)-li(x)|=O(x^{0.5+\varepsilon}),\quad \varepsilon>0$ of Riemann Hypothesis is fully established. Riemann Hypothesis is proved!!**

## 4. Prove another equivalent proposition of RH $|\psi(x)-x|=O(x^{0.5+\varepsilon}),\quad \varepsilon>0$

We make more transformation to $\zeta(s)$: $\dfrac{1}{s}M_e(s)=-\dfrac{1}{s}\zeta'(s)=-\dfrac{1}{s}(1+\dfrac{1}{2^s}+\dfrac{1}{3^s}+\dfrac{1}{4^s}+...)'=$ $-\dfrac{1}{s}(1+e^{-s\log 2}+e^{-s\log 3}+e^{-s\log 4}+...)'=\dfrac{1}{s}(\dfrac{\log 2}{2^s}+\dfrac{\log 3}{3^s}+\dfrac{\log 4}{4^s}+...)=L(M_e(x))$. The real function corresponding to this complex function is the logarithmic contraction from abscissa of the function as: $M(x)=\sum_{m=1,2,..}^{m\le x}\log m$, so the new function can be recorded as: $M_e(x)$. In number theory, there is an important function: $\psi(x)=\sum_p u(x-p^k)\cdot\log p,\quad p=prime\ number$; and a famous formula:

$\sum_{k=1,2,3,..}^{k\le x}\psi(x/k)=\sum_{m=1,2,...}^{m\le x}\log m$ **...... (4.1)** It can be corresponded to complex functions relation:

By Euler product formula, $\zeta(s)=\Pi(1-p^{-s})^{-1}$, we let $\psi_e(x)=\psi(e^x),\quad \psi_e(\log x)=\psi(x)$.

$$L[\psi_e(x)]=\frac{1}{s}\Psi_e(s)=-\frac{1}{s}[\log\zeta(s)]'=-\frac{1}{s}[\log\Pi(1-p^{-s})^{-1}]'=-\frac{1}{s}\sum(p^{-s}+\frac{1}{2}p^{-2s}+\frac{1}{3}p^{-3s}+...)'=$$

$-\dfrac{1}{s}\sum(e^{-s\log p}+\dfrac{1}{2}e^{-2s\log p}+\dfrac{1}{3}e^{-3s\log p}+...)'=\dfrac{1}{s}\sum(\dfrac{1}{p^s}+\dfrac{1}{p^{2s}}+\dfrac{1}{p^{3s}}+...)\cdot\log p$ . On the other hand,

$$-\frac{1}{s}[\log\zeta(s)]' = -\frac{\zeta'(s)}{s\zeta(s)}, \quad \text{so:} \quad \frac{\Psi_e(s)}{s} = \frac{M_e(s)}{s\zeta(s)}, \quad \frac{1}{s}\zeta(s)\cdot\Psi_e(s) = \frac{1}{s}M_e(s) \quad \text{......} \; \mathbf{(4.2)}$$

**Figure 5**

When the abscissa coordinates logarithms, $\log(x/k) = \log x - \log k$ is essentially a convolution sampling sum. The sampling function and the sampled function can be exchanged: When $\psi_e(x)$ is considered to be a sampled function, the sampling pulse sequence is $zeta(x)$; if $\int_0^x zeta(t)dt$ is regarded as a sampled function, the sampling pulse sequence is $L^{-1}[\Psi_e(s)]$. Figure 5 is the first case.

Using the relational formula **(4.1)** above, suppose $\psi(x) = x + A\log x + h(x)$, so $\sum_{k=1,2,...}^{k\le x}\psi(x/k) =:$

$$\sum_{k=1,2,...}^{k\le x} x/k + A\sum_{k=1,2,...}^{k\le x}\log(x/k) + A\sum_{k=1,2,...}^{k\le x}h(x/k) = x\log x + \gamma\cdot x + \frac{1}{2} + A\left(x\log x - \sum_{k=1,2,...}^{k\le x}\log k\right) + A\sum_{k=1,2,...}^{k\le x}h(x/k)$$

$$= x\log x + \gamma\cdot x + \frac{1}{2} + A(x\log x - x\log x + x - \frac{\log x}{2} - \frac{\log 2\pi}{2}) + A\sum_{k=1,2,...}^{k\le x}h(x/k) = \sum_{m=1,2,...}^{m\le x}\log m = x\log x - x + \frac{\log x}{2} + \frac{\log 2\pi}{2}$$

, get $\gamma + A = -1$, $A = -1-\gamma$. Obviously, the mean function of $\psi(x)$ is $x-(1+\gamma)\log x$, and its vibration term $h(x)$ is sampled and sum up to $\sum_{k=1,2,...}^{k\le x}h(x/k) \approx \frac{1}{1+\gamma}\left(\frac{\gamma\log x}{2} + \frac{\gamma\log 2\pi}{2} + \frac{1}{2}\right)$.

The mean function of $\psi(x)$ can no longer increase or decrease any constant c, because in that case, $\sum\psi(x/k)$ will increase or decrease $cx$, which is totally wrong. As long as the mean function is contracted by abscissa in logarithm, the mean value function of $\psi_e(x)$ is obtained: $e^x-(1+\gamma)x$.

Further, let: $\int_0^x \psi_{e0}(x)dx = e^x - (1+\gamma)x$, $\psi_{e0}(x) = e^x - (1+\gamma)$ **...... (4.3)** , called the ***original function of*** $\psi_e(x)$ ***`s mean function.***

The original function of $\psi_e(x)$ is: $\psi'_e(x) = \sum \log p \cdot \delta(x - k\log p)$, the original function of $J_e(x)$ is: $J'_e(x) = \sum (1/k)\cdot\delta(x - k\log p)$. There is a relationship between them:

$\psi'_e(x) = xJ'_e(x)$. If $\psi'_e(x)$ and $J'_e(x)$ are averaging, such relationship still exists, so there is:

$xJ_{e0}(x) = \psi_{e0}(x) = e^x - (1+\gamma),\quad J_{e0}(x) = \dfrac{e^x}{x} - \dfrac{1+\gamma}{x}$ **...... (4.4)** We can check whether the mean function is valid. The result can be obtained by the integration by parts as:

$\int \psi'_e(x)dx = \int xJ'_e(x)dx = \dfrac{x^2}{2}J'_e(x) - \int \dfrac{x^2}{2}J''_e(x)dx$, use $J_{e0}(x)$ to replace $J'_e(x)$,

$$\frac{x^2}{2}\left(\frac{e^x}{x} - \frac{1+\gamma}{x}\right) - \int \frac{x^2}{2}\left(\frac{e^x}{x} - \frac{1+\gamma}{x}\right)' dx = \frac{xe^x}{2} - \frac{(1+\gamma)x}{2} - \int \frac{x^2}{2}\left(\frac{xe^x - e^x}{x^2} + \frac{1+\gamma}{x^2}\right)dx = e^x - (1+\gamma)x.$$

So, $J_{e0}(x) = \dfrac{e^x}{x} - \dfrac{1+\gamma}{x}$ is surely the original function of $J_e(x)$ `s mean function.

From $lie(x) = \int_1^x \dfrac{e^t}{t}dt$, we get $lie(x) - \int_1^x J_{e0}(t)dt = \int_1^x \dfrac{1+\gamma}{t}dt = (1+\gamma)\log x$, obviously the order is very low, indicating that the two functions are very close; $lie(x)$ is slightly larger than the mean value of $Je(x)$, so $Je(x)$ 's oscillation around $lie(x)$ will be more times under $lie(x)$.

***Another important equivalent proposition of RH*** $|\psi(x) - x| = O(x^{0.5+\varepsilon_1}),\quad \varepsilon_1 > 0$ ***has been proved.*** **It requires:** $|\psi_e(x) - e^x| = O(e^{0.5x+\varepsilon_2}),\quad \varepsilon_2 > 0$. Because the original function of $\psi_e(x)$ `s mean function has been determined by (4.3) formula. ***The functions of the upper and lower bounds of*** $\psi_e(x)$ ***are also some kinds of mean function. So, the original functions of the upper and lower bounds of*** $\psi_e(x)$ ***and*** $Je(x)$ ***must have the same relations as their mean value functions, that means:***

$\psi^+_{e0}(x) = xJ^+_{e0}(x)$ and $\psi^-_{e0}(x) = xJ^-_{e0}(x)$. We have got $|Je(x) - lie(x)| < \dfrac{x}{2} + \dfrac{(x/2)^2}{2\cdot 2!} + \dfrac{(x/2)^3}{3\cdot 3!} + \ldots < \int_1^x \dfrac{e^{0.5t}}{t}dt$ in part 3 and $lie(x) = \int_1^x \dfrac{e^t}{t}dt$ in part 2, so there should be:

$$-\int_1^x \frac{e^{0.5t}}{t}dt < \int_1^x \frac{\psi_e'(t)}{t}dt - \int_1^x \frac{e^t}{t}dt < \int_1^x \frac{e^{0.5t}}{t}dt, \quad \int_1^x \frac{e^t - e^{0.5t}}{t}dt < \int_1^x \frac{\psi_e'(t)}{t}dt < \int_1^x \frac{e^t + e^{0.5t}}{t}dt,$$

$$\int_1^x (e^t - e^{0.5t})dt < \int_1^x \psi_e'(t)dt < \int_1^x (e^t + e^{0.5t})dt, \quad e^x - 2e^{0.5x} < \psi_e(x) < e^x + 2e^{0.5x},$$

$$\left|\psi_e(x) - e^x\right| < 2e^{0.5x}, \text{ that means } \left|\psi(x) - x\right| < 2\sqrt{x} = O(x^{0.5+\varepsilon}), \quad \varepsilon > 0$$

The process of the proof in part 3 can show ***why the vibration range of*** $Je(x) - lie(x)$ ***is at the*** $O(e^{0.5x})$ ***order*** too: The fundamental reason is that the exact upper and lower bounds of $zeta_1(x)$ are $e^x$ and $e^x - 1$, that means, the average error of $r(x)$ is about 0.5. The 1/2 error constitutes the $1/2s$ term of the ***convolution equivalent function*** (see part 3), and is eventually transformed into the exponential function amplitude of $x/2$, that is, the order $O(e^{0.5x})$.

*Why many "Generalized Riemann Hypothesis" also assume that the real parts of "non trivial 0 points" are all +1/2? Because* ***the functions in these conjectures can be considered as a combination of some arithmetic transformations of the Zeta function.*** *These transformations can change the size of the main item and the period and phase of the vibration, but* ***it does not change the law that the maximum difference between the function itself and its mean function is +/-0.5.*** *For two examples:*

$$\frac{Ze_2(s)}{s} = \frac{1}{s}\left(\frac{1}{2^s} + \frac{1}{4^s} + \frac{1}{6^s} + \frac{1}{8^s} + \ldots\right), \qquad \frac{Ze_{1.5}(s)}{s} = \frac{1}{s}\left(\frac{1}{1.5^s} + \frac{1}{2.5^s} + \frac{1}{3.5^s} + \frac{1}{4.5^s} + \ldots\right)$$

Calculations show when $2 \le x \le 2 \cdot 10^{14}$, $\left|J(x) - Li(x)\right| < \frac{0.7\sqrt{x}}{\log x}$, as in Figure 6.

Michael Atiyah <M.Atiyah@ed.ac.uk>
周二 2018/9/25 22:19
你

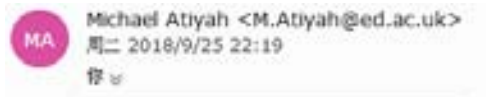

zhu min jing,
The zeta function has a pole at s=1, it is analytic on the plane punctured at s = 1. But this punctured plane is not convex so zeta is polynomial on compact convex subsets of the critical strip. The same is true for Todd and the composition of Todd and zeta.

There are many weakly analytic functions, but Hirzebruch showed that, among them, Todd could be uniquely characterized by a simple property. Read chapter 1 of his book for the full story.

Your work on RH has some overlaps with my simpler version, but I leave it to you to examine the comparison,

Michael Atiyah
On 25/09/2018 12:53, Zhu Jing Min wrote:

> Dear Sir,
>
> Excuse me, I am reading your proof of RH and have some questions to ask.
>
> 1.Do you think that both the Zeta function and the Todd function are belong to weak analytic functions, while the Todd function is more general than the Zeta function, and the Zeta function belongs to the set of Todd function, or is it a type of Todd function? So Zeta functions must have the properties of Todd functions?
>
> 2.If the answer to these questions is "yes", what is your basis?
>
> 3.You said in the previous emails that "I have now had time to look at your paper on the Riemann Hypothesis. I recognize in there many features that appear in my proof of RH", which features do you mean in my paper? And how they can appear in your proof?
>
> Thanks & Regards
>
> Zhu Jing Min
>
> ---
>
> **发件人:** Michael Atiyah <M.Atiyah@ed.ac.uk>
> **发送时间:** 2018年9月24日 23:35
> **收件人:** Zhu Jing Min
> **主题:** Re: Could you kindly give me a help?
>
> zhu jing min,
>
> I regret that shortage of time meant that I omitted to mention your name. I will rectify that in the published version in my own words. I do not use social media.
>
> best wishes,
>
> Michael Atiyah
>
> On 24/09/2018 16:08, Zhu Jing Min wrote:
>
>> **Dear Sir Michael F Atiyah,**
>>
>> I have watched your lecture. It is wonderful and inspiring. Anyway, I have not been mentioned in it. Could you kindly help me by writing a few words on Twitter or Facebook, such as (of course, any other words you think them Ok):
>>
>> A man from China have offered his proof of Riemann Hypothesis too. He say there are some new ideas in his paper. It can be reached by:
>>
>> https://arxiv.org/abs/1807.00849
>>
>> Thank you very much!
>>
>> Best Regards
>>
>> Zhu Jing Min
>>
>> ---
>>
>> **发件人:** Zhu Jing Min <jing_min_zhu@hotmail.com>
>> **发送时间:** 2018年9月23日 21:53
>> **收件人:** Michael Atiyah
>> **主题:** Re: Could you kindly give me a help?
>>
>> **Dear Sir Michael F Atiyah,**
>>
>> I am so glad to hear from you! And you can give me your proof, I am very excited! Thank you for your help again!! I will seriously study your proof, hoping to understand it. Of course, I will not disclose your proof until you say it can be opened.
>>
>> I think there should be some different solutions to a problem, especially for Riemann Hypothesis, which is well known for its many equivalence propositions.
>>
>> According to Chinese grammatical order, Zhu is my family name and Jing Min is given name.
>>
>> When you mention my name in the speech tomorrow, please give the web address of my article on arXiv. That is:
>>
>> https://arxiv.org/abs/1807.00849
>>
>> I hope more mathematicians will pay attention to my approach and discuss it in depth, which I think will be useful for the research of number theory.
>>
>> Of course, my proof may turn out to be wrong in the end, but it doesn't matter. I'm not ashamed while so many people were wrong in the past. But I really don't want my job to be buried.
>>
>> Finally, I sincerely wish you a successful speech tomorrow!
>>
>> Sincerely
>>
>> Zhu Jing Min
>>
>> ---
>>
>> **发件人:** Michael Atiyah <M.Atiyah@ed.ac.uk>
>> **发送时间:** 2018年9月23日 18:36
>> **收件人:** Zhu Jing Min
>> **主题:** Re: Could you kindly give me a help?
>>
>> Zhu Jing Min,
>> I have now had time to look at your paper on the Riemann Hypothesis. I recognize in there many features that appear in my proof of RH, but your approach mixes computation and proof. My approach is much simpler.
>> You go back to real numbers, I go forward to quaternions.
>> However I do not directly attack computation, but leave that to another stage.
>> Your ideas may have some use at this stage.
>> I attach in a subsequent email a copy of my paper, which I will present tomorrow here in Heidelberg. I rely on you to keep it confidential until then. In return I will mention your name in my lecture.
>> Best wishes,
>> Michael Atiyah
>>
>> On 22/09/2018 19:47, Zhu Jing Min wrote:
>>
>>> **Dear Sir Michael F Atiyah:**
>>>
>>> Excuse me, this is the third time I wrote to you. I also promise this is the last time. I come from Guangzhou, China, and have studied Riemann Hypothesis for 10 years. I think I have found a completely new way, which has never been tried by others before, and obtained some important results. But because I am not a well-known man in mathematics society, so far my papers have not been taken seriously, which makes me very anxious and distressed. I am worried that my research will be buried, and I believe that even if my paper is not 100% correct, there are still some new discoveries worth discussing and studying.
>>>
>>> Could you spend a little time reading this paper? I guarantee it is very easy to understand. If you think there are some innovations and discoveries too, could you write a litte bit about my work on your Twitter or Facebook accounts? Or would you like to mention a little bit about my work in your upcoming speech on Sep. 24? If so, I will be very grateful and thankful to you!!
>>>
>>> The pdf file of my paper is attached. It also can be reached at:
>>>
>>> https://arxiv.org/abs/1807.00849
>>>
>>> Many new ideas, only 14 pages(core part<10 pages), 7 figures.
>>>
>>> **Summary:**
>>>
>>> In the past 100 years, the research of Riemann Hypothesis meets many difficulties. Such situation may be caused by that people used to study Zeta function only regarding it as a complex function. Generally, complex functions are far more complex than real functions, and are hard to graph. So, people cannot grasp the nature of them easily. Therefore, it may be a promising way to try to correspond Zeta function to real function so that we can return to the real domain to study RH.
>>>
>>> In fact, by the view of Laplace transform, the whole picture of Zeta(s)/s is very clear and simple(very closing exp(x)), and the problem can be greatly simplified. And by Laplace transform, most integral and convolution operations can be converted into algebraic operations, which greatly simplifies calculating and analysis.
>>>
>>> Best Regards
>>>
>>> Zhu Jing Min

**Figure 6**

$\dfrac{J(x)-Li(x)}{\sqrt{x}/\log x}$

+0.4

−0.4

1 4 7 10 13 16 19 22 25 28 31 34 37 40 43 46 49 52 55 58 61 64 67 70 73 76 79 82 85 88 91 94

from 2 to 2*10^14, step 0.5log2

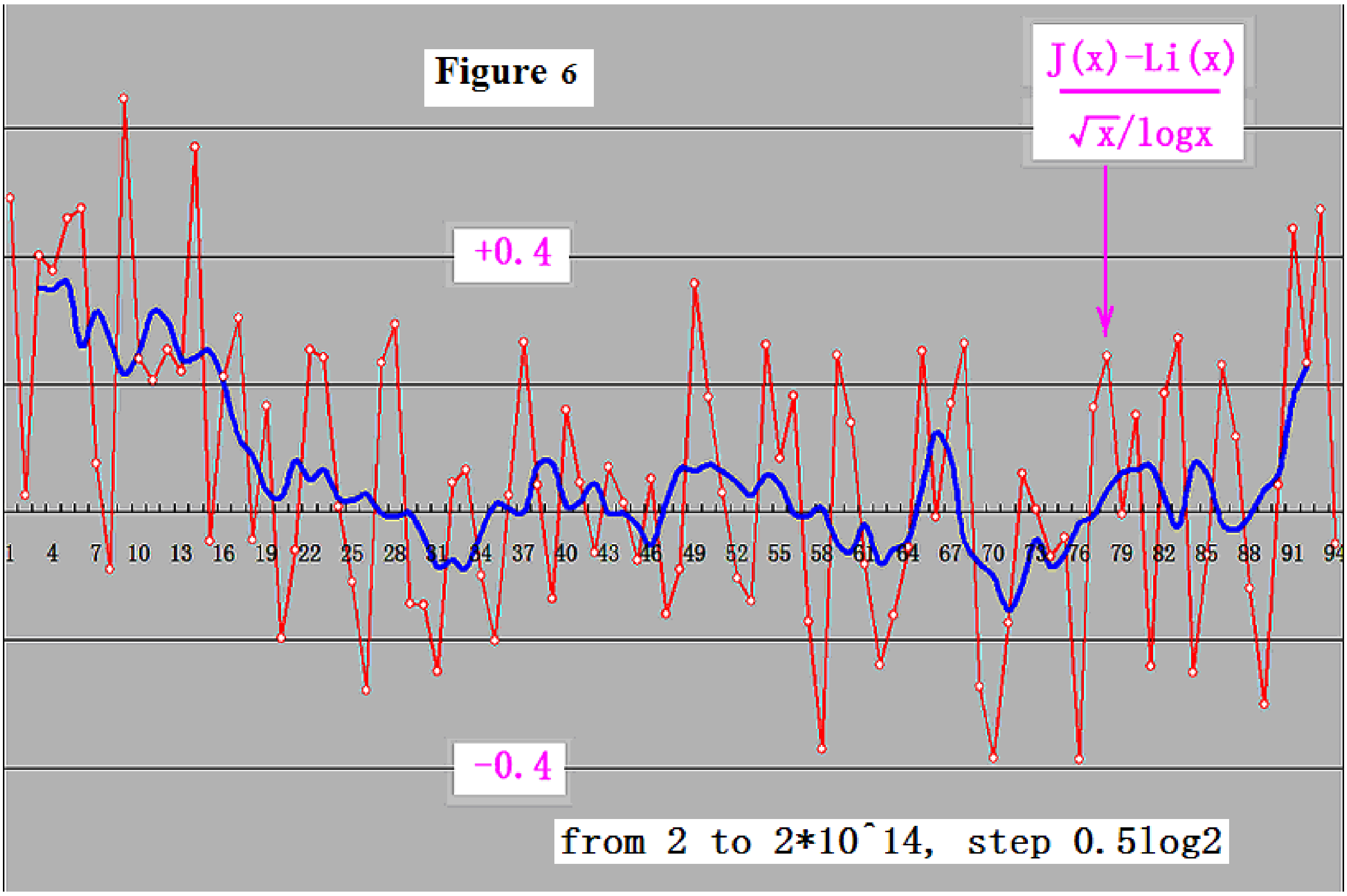

## A brief introduction of the author

I was born in Guangzhou, China in December 1972. I studied at South China University of Technology, major in radio engineering. After graduation, I had worked in a few of multinationals and Chinese companies, in charge of technique, sales, and management, and got MBA degree in 2010. I am interested in science and technology, good math's skills, broad thinking and good at summing up experience.

In 2008, I began to research the distribution of prime numbers. In 2009, began to study Riemann Hypothesis. From 2010, the Zeta function and other related problems were investigated from the view of Laplace transformation. There are new and important discoveries and many efforts to overcome the difficulties. Finally, in April 2018, I find a way to solve the RH and complete the first proof.

**Contact: jing_min_zhu@hotmail.com**

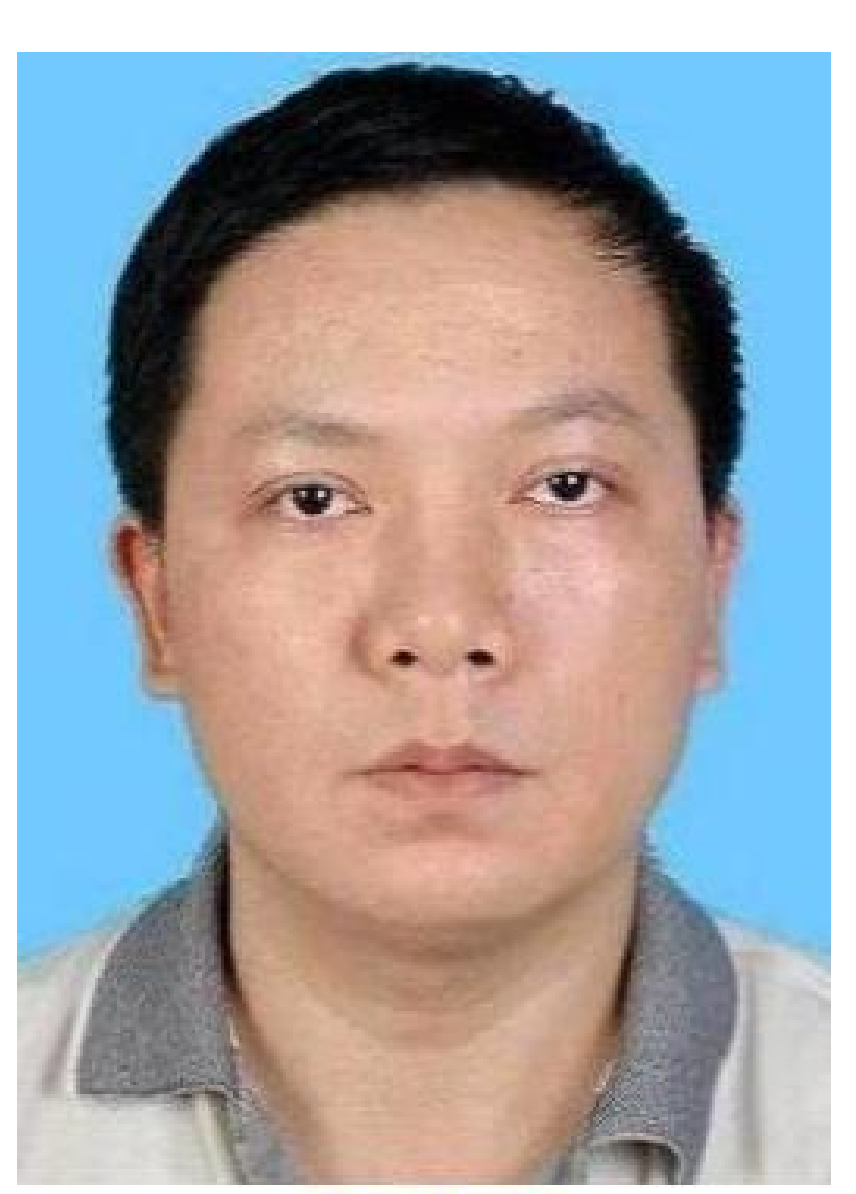

朱敬民